\documentclass{ifacconf}
\usepackage{amsmath}
\usepackage{amssymb}
\usepackage{breqn}

\usepackage{bm}
\usepackage{algorithm}
\usepackage[]{algorithmicx}
\usepackage[]{algpseudocode}
\usepackage{setspace}
\usepackage{varwidth}

\usepackage{float} 

\usepackage{adjustbox}

\usepackage[intoc]{nomencl}
\makenomenclature

\usepackage{enumerate}
\usepackage{graphicx}      
\usepackage{natbib}        
\begin{document}
\begin{frontmatter}

\title{Data Assimilation for Combined Parameter and State Estimation in Stochastic Continuous-Discrete Nonlinear Systems\thanksref{footnoteinfo}} 

\thanks[footnoteinfo]{Authors T. Diaa-Eldeen, C. F. Berg, and M. Hovd are with BRU21 – NTNU Research and Innovation Program on Digital and Automation Solutions for the Oil and Gas Industry (www.ntnu.edu/bru21) and are supported by Equinor.}

\author[NTNU-ITK]{Tarek Diaa-Eldeen} 
\author[DTU]{Marcus Krogh Nielsen}
\author[NTNU-IGP]{Carl Fredrik Berg}
\author[NTNU-ITK]{Morten Hovd}
\author[DTU]{John Bagterp Jørgensen}

\address[NTNU-ITK]{Department of Engineering Cybernetics, Norwegian University of Science and Technology (NTNU), N-7491 Trondheim, Norway}
\address[NTNU-IGP]{Department of Geoscience and Petroleum, Norwegian University of Science and Technology (NTNU), N-7491 Trondheim, Norway}
\address[DTU]{Department of Applied Mathematics and Computer Science, Technical University of Denmark, DK-2800 Kgs Lyngby, Denmark }

\begin{abstract} Data assimilation (DA) provides a general framework for estimation in dynamical systems based on the concepts of Bayesian inference. This constitutes a common basis for the different linear and nonlinear filtering and smoothing techniques which gives a better understanding of the characteristics and limitations of each approach. In this study, four nonlinear filters for systems governed by stochastic continuous dynamics with discrete measurements are inferred as different approximate solutions to the Fokker-Planck equation in the prediction step of the Bayesian update. The characteristics and drawbacks of each filtering approach are discussed in light of the introduced approximations to the original Bayesian problem in each case. The introduced continuous-discrete (CD) filters are then implemented to solve the combined estimation problem of the state and parameter variables in a nonlinear system described by a stochastic differential equation (SDE), namely, the continuous stirred-tank reactor (CSTR). The performance measures for each filter, such as estimation errors, computational cost, and the ability to quantify the correct uncertainty, are compared and investigated with attribution to the assumed approximations introduced in the derivation of each class of filters.
\end{abstract}


\end{frontmatter}

\section{Introduction}
Data assimilation (DA) is the theory of combining prior knowledge from dynamical systems with new observational information in an optimal way to obtain the best description of the underlying dynamical system, and characterize its uncertainty \citep{carrassi2018data}. Devised based on Bayesian inference, DA is considered a general estimation theory that introduces different methods to optimally fuse all the available information to update a prior uncertainty distribution. Although the term data assimilation was first used in numerical weather predictions and geosciences, its mathematical foundation is rooted in Bayesian inference, control theory, and variational calculus \citep{evensen2022data}. Two main advantages arise from formulating the estimation problem this way. First, it provides a thorough framework where the different estimation methods can be derived and formulated from the same origin. This top-down approach, which is based on Bayes theorem, can be applied to different estimation problems involving different uncertain variables. For example, estimating uncertain state variables, parameters, control inputs (forcing inputs), and disturbances, in prediction, filtering, or smoothing problems. Secondly, this Bayesian formulation allows to avoid solving the estimation problem by directly inverting the observations, which can be problematic in ill-posed and high-dimensional problems \citep{Diaa-Eldeen}; instead, it uses the likelihood of the observations to update the prior uncertainties of the unknown variables.

From an optimization-based point-of-view, Bayes' formula is optimal in the sense that it constitutes the optimal solution that balances the different sources of information in the problem. In addition to its direct derivation from the definition of the joint probability between the uncertain variables and the observations, Bayes' theorem can be derived from solving three different optimization problems using those sources of information
\begin{enumerate}[(i)]
    \item The joint minimization of the relative entropy between the posterior and the prior, represented by the Kullback-Leibler (KL) divergence, and the mean squared errors (MSE) from the likelihood of the observed data. The first-order optimality condition of this optimization problem constitutes the Bayes’ formula, with a unique solution that precisely describes the posterior distribution \citep{BayesOptimal}.
    \item Minimization of the output and input information discrepancy \citep{OptimalInformation}.
    \item Minimization of the upper bound of a variational inference problem \citep{bui2021optimality}.
\end{enumerate}
Solving these optimization problems results in the same well-known Bayes formula
\begin{equation} \label{eq:Bayes}
p(X|Y) = \frac{p(Y|X) p(X)} {p(Y)},
\end{equation}
where the posterior $p(X|Y)$ is the solution to the assimilation problem, $p(Y|X)$ is the likelihood function, $p(X)$ is the prior, and $p(Y)$ is the normalization term.

The appealing feature in DA is that using different assumptions and approximations, all the various estimation and filtering methods, to name a few, the Kalman Filter (KF), Extended KF (EKF), Unscented KF (UKF), Ensemble KF (EnKF), and the Particle Filter (PF), can be derived from the same starting point\textemdash Bayes theorem. This helps to discriminate the characteristics and the limitations of the different methods in a neat and comprehensive manner.

In this paper, four continuous-discrete (CD) nonlinear filters that are widely applied in control systems literature, namely, the CD-EKF, CD-UKF, CD-EnKF, and CD-PF, are formulated using the Bayesian DA framework. The inherent performance characteristics and limitations of each filtering approach, which mainly emerge from the imposed assumptions and approximations to the original Bayesian problem, are discussed and empirically emphasized and investigated in a numerical case study. The different CD nonlinear filters are implemented to update both the state and parameter variables in a stochastic nonlinear model of a chemical process involving an exothermic reaction in an adiabatic continuous stirred tank reactor (CSTR).

The following parts of this paper are organized as follows. Section \ref{sec:BayesCS} introduces Bayesian inference and estimation in CD systems. In Section \ref{sec:NLF}, the different nonlinear filtering approaches are introduced as approximate solutions to the Bayesian problem. Then, Section \ref{sec:App} shows a numerical application of the filtering approaches to a CSTR system modeled with stochastic continuous nonlinear dynamics. Finally, Section \ref{sec:conc} concludes the paper.

\section{Bayesian Inference in Nonlinear Continuous-Discrete Systems} \label{sec:BayesCS}
\subsection{Nonlinear Continuous-Discrete Stochastic Models}
Consider a continuous-time nonlinear dynamical system that can be modeled by a discretely-observed stochastic differential equation (SDE), with the stochastic term defined in the sense of Itô integral 
\begin{subequations}
\begin{equation}
    \begin{split}
    dx(t) &= f(t,x(t),u(t),d(t),\theta)dt \\
    &+ \sigma(t,x(t),u(t),d(t),\theta)d\omega(t),
    \label{eq:StochModel}
    \end{split}
\end{equation}
\begin{equation}
    y(t_k) = h(t_k,x(t_k),\theta) + v(t_k),
    \label{eq:ObsModel}
\end{equation}
\label{eq:CDStochModel}
\end{subequations}
where $f(\cdot)$ is the nonlinear model operator representing the drift term, $\sigma (\cdot)$ is the diffusion function such that the diffusion tensor $D = {1\over2} \sigma \sigma^T$, and $h(\cdot)$ is the nonlinear measurement operator. $x(t) \in \mathbb{R}^{n_x}$, $u(t) \in \mathbb{R}^{n_u}$, $d(t) \in \mathbb{R}^{n_d}$, and $\theta \in \mathbb{R}^{n_\theta}$ represent the model state, input, disturbance, and parameter vectors, respectively, as functions of time. The process noise $\omega(t) \in \mathbb{R}^{n_\omega}$ is a standard
Wiener process with a unity diffusion matrix, $Q_c(t) = I$, such that $d\omega \sim \mathcal{N}_{iid}(0, Idt)$. $v(t_k) \sim \mathcal{N}_{iid}(0, R_k)$ represents the measurement noise term, which includes instrumentation errors and also representation errors that can result in from unresolved processes in the available models. 

For the sake of the joint estimation of different unknown variables, an augmented state-vector representation is used. In this study, combined state and parameter estimation is considered, and the augmented state vector, $\hat{x} \in \mathbb{R}^{n_{\hat{x}}}$, is defined as $\hat{x}(t) = \left(\begin{smallmatrix}x(t) \\ \theta(t)\end{smallmatrix}\right)$, where the parameter variable, $\theta(t)$, is described by \eqref{eq:CDStochModel} with $f(.) = 0$. In addition, only additive stochasticity is considered, with a state-independent diffusion term.

\subsection{Classification of Data Assimilation Methods}
DA methods can be classified into two main categories depending on the estimation procedure which, can be implemented either in batches (variational methods) or recursively (sequential methods).

\subsubsection{Variational Methods.} The estimation problem is cast as an optimization problem that is solved over the assimilation horizon using all the available observations. An example of this is the Moving Horizon Estimation (MHE) method. The variational (batch) estimation methods solve the original Bayesian problem in \eqref{eq:Bayes}, where the joint distribution of all the states is conditioned on all the available observations. However, two main drawbacks arise in this approach---the computation time is inefficient for real-time applications, and they do not provide a means to quantify the uncertainty in the estimates.

\subsubsection{Sequential Methods.} \sloppy Introducing two assumptions to the original Bayesian problem allows to implement it recursively in time. That is, assuming that the dynamical model \eqref{eq:StochModel} is a first-order Markov process, and the measurement models \eqref{eq:ObsModel} are conditionally independent (with uncorrelated measurement errors in time), \emph{i.e.}, $p(Y|X) = \Pi_k p(y_k|x(t_k))$, where $y_{1:k} = \{y_1, y_2, ..., y_k\}$, allows to solve the Bayesian problem in a recursive manner when new information is available. In other words, $p(x(t_k)|y_{1:k}) \sim p(x(t_k))\Pi_k p(y_k|x(t_k))$, which is a practical updating approach in real-time applications. However, this recursive update can give two different types of solutions. The original solution that is obtained from the introduced assumptions is the smoothing solution, where the joint distribution of the state variables is recursively conditioned on all the observations. That is, the posterior distribution is defined as $p(X|Y) \overset{\Delta}{=} p(x(t)|y_{1:k})$, for $t_0 < t \le t_k$. This means that the new observational information at time $t_k$ updates the current state, $x(t_k)$, and it is also projected backward in time to update the state variables from time zero. Alternatively, marginalizing the joint distribution to the current state at observation time $t_k$ by integrating out the state trajectory $x(t)$ except at the current assimilation step gives the filtering solution $p(x(t_k)|y_{1:k})$. Both the smoothing and the filtering solutions are equivalent at the end of the assimilation horizon, $x(t_k)$. The smoother in addition updates the state trajectory at previous times. On the other hand, for the same assimilation window, the smoothing solution is equivalent to the variational solution in the linear Gaussian case.

Using this continuous-state discrete-measurements formulation makes it possible to obtain the Bayesian posterior at any time when new observational information is available, not necessarily at fixed assimilation steps, which is a required property in many practical applications. The next section discusses the uncertainty propagation to the next observational time in CD systems.

\subsection{Uncertainty Propagation and Error Evolution}
In spite of the ostensible simplicity of the Bayesian formula, a closed-form solution can rarely be found. This is regarded to the complexity of the required computations to solve the involved integrals. For example, in CD systems, the Bayesian update can be implemented in two or three steps. The first step is the prediction step (time update), where a prior distribution, $p(x(t)|y_{1:k-1})$ for $t>t_{k-1}$, is calculated by propagating an initial probability distribution, $p(x(t_{k-1})|y_{1:k-1})$, forward in time. The exact solution for this propagation is determined by solving a boundary-value problem given by the Fokker-Planck equation (FPE) (Kolmogorov's forward equation) \eqref{eq:FP} with the initial probability distribution representing the boundary condition. This can be the initial uncertainty that characterizes our belief in the unknown variables, $p(x(t_0))$, or the posterior from the last assimilation step, $p(x(t_{k-1})|y_{1:k-1})$. The solution of the FPE gives the time-evolution of the probability density function (pdf) for a general stochastic model \eqref{eq:CDStochModel}. Then, the measurement update step, where this prediction is conditioned on the new observational information at time $t_k$ represented by the likelihood function, $p(y_k|x(t_k))$, to update the uncertainty and solve for a posterior distribution, $p(x(t_k)|y_{1:k})$.
The FPE is
\begin{equation}
    \partial_t p = -\nabla \cdot J(t,\hat x(t),u(t),d(t),\theta),
    \label{eq:FP}
\end{equation}
where the probability current,  $J$, for \eqref{eq:StochModel} is defined as \citep{Thygesen:2022,Pavliotis:2014}
\begin{equation}
    J = f p - \nabla \cdot (D p). 
\end{equation}
Alternatively, the time update step in CD systems can be implemented in two separate steps---a discretization step, where the transition density, $p(x(t_k)|x(t_{k-1}))$, is solved for by integrating the FPE, and an update step, where the Chapman-Kolmogorov equation is applied to compute the predicted probability, $p(x(t_k)|y_{1:k-1})$, \citep{SimoSarrkaPhD}.

Due to the complexity of the resulting integrations that tend to be in high-dimensional spaces, an exact solution is usually not feasible or avoided. Therefore, the different filtering approaches can be derived as equivalent or approximate solutions to the FPE such that each filter will belong to one (or more) of the following cases.

\subsubsection{Exact Equivalent Solution.} An exact equivalent solution to the FPE can be found in the case of linear systems with Gaussian noises, where the probability distribution is completely determined by the first two statistical moments. In this case, the exact prior can be computed by propagating the mean and the covariance information using the linear dynamics, resulting in the standard KF. This is equivalent to solving the FPE for the first two statistical moments. Otherwise, further numerical approximations are required to find an approximate solution.

\subsubsection{Linearization and Gaussian Approximation.}
In the nonlinear case, the Bayesian posterior can be approximated with a Gaussian distribution, $p(x(t)|y_{1:k}) \approx \mathcal{N}(x(t)|m(t), P(t))$, where $m(t)$ and $P(t)$ are the mean and the covariance of the approximated Gaussian process, resulting in the Kalman-based nonlinear filters. This can be done either using moment matching or by minimizing the discrepancy between the exact and approximate processes. A linearized model can be used for this purpose to propagate the Gaussian process information. The linearization can either be dynamical linearization, which results in the EKF, or derivative-free statistical linearization, which is the case in the statistically-linearized filters. In this case, the propagated distribution is guaranteed to be Gaussian, given that the initial model errors are Gaussian. However, the propagation is implemented using a linear function approximation. Alternatively, a Gaussian approximation for the posterior can be obtained by explicitly approximating the updated uncertainty by a Gaussian distribution, which is the case in other derivative-free nonlinear Kalman-based filters, such as the EnKF, the UKF, and other sigma-point filters. In this case, the original nonlinear model is used to propagate the error statistics, and a Gaussian process is approximated after propagation.

\subsubsection{Sampling Approximation.} Sampling is an efficient alternative to propagate the uncertainty and approximate the solution to the FPE. That is, instead of solving for the exact distribution of the propagated error statistics, it can be approximated by propagating a finite number of model realizations drawn from an initial distribution, $x^{(i)}(t_{k-1}) \sim p(x(t_{k-1})|y_{1:k-1})$. This can be implemented either deterministically or using Monte Carlo (MC) random sampling. The latter is equivalent to using a Monte Carlo Markov Chain (MCMC) approach to integrate the FPE. The deterministic sampling approach is the error propagation approach used in the UKF and other sigma-point filters, while the MC-based sampling approach is used in the EnKF and PF. In the UKF and the EnKF, a Gaussian process approximation in addition is applied after the propagation such that the pdf is described by the sample mean and the sample covariance. While in the PF, the full probability distribution is sampled.

Other approaches for approximating the solution of the FPE include using the Galerkin method as in the nonlinear projection filters \citep{GalerkinMethod}, and deep learning \citep{xu2020solving}. 

\section{Nonlinear Filters as Data Assimilation Approaches} \label{sec:NLF}
Given an initial distribution, $p(x(t_0))$, a measurement sequence, $y_{1:k} = \{y_1, y_2, ..., y_k\}$, and a CD probabilistic state space model which describes the dynamic model, $x(t) \sim p(x(t)|x(t_0))$, for $t_0 < t \le t_k$, and the measurement model, $y_k \sim p(y_k|x(t_k)) \overset{\Delta}{=} \mathcal{N}(y_k|h(x(t_k),t_k),R_k)$, as given by \eqref{eq:StochModel} and \eqref{eq:ObsModel}, respectively, the different filters can be implemented using the previously introduced approximations to sequentially solve the FPE and update the Bayesian posterior.

First, The CD-KF constitutes the closed-form solution to the Bayesian problem in the linear Gaussian case. The solution is found directly by applying the two-step approach---prediction, using the exact equivalent solution to the FPE, and then the update step to compute the Bayesian posterior. Alternatively, the three-step approach can also be applied in this case by adding a discretization step, $p(x(t_k)|x(t_{k-1})) = \mathcal{N}(x(t_k)|A_{k-1}x(t_{k-1}), Q_{k-1})$, where $A$ and $Q$ are the discretization of the system matrix and the excitation noise covariance matrix of the linear system, respectively; and then implementing the time and measurement updates similar to the discrete-time KF. Otherwise, the introduced numerical approximations are applied as follows.

\subsection{The Extended Kalman Filter}
In the nonlinear-based Kalman filters, the moment equations for the first and second moments are solved approximately. The CD-EKF belongs to the first class of approximate solutions to the FPE, where a linearized model step is used to approximate the propagation of the pdf by a Gaussian process that has the same mean and covariance as the original process at the linearization point. Taylor-series is used to obtain the linear or the quadratic approximation for the nonlinear SDE in \eqref{eq:CDStochModel}, resulting in the first- or the second-order CD-EKF, respectively. The time and measurement update for the first-order CD-EKF are then implemented as shown in Algorithm \ref{Alg:CD-EKF}, where the Jacobain matrices $A_k(t) = \frac{\partial f(t)}{\partial \hat{x}}|_{\hat{x} = \hat{x}(t_k)}$ and $C_k = \frac{\partial y(t_k)}{\partial \hat{x}}|_{\hat{x} = \hat{x}(t_k)}$ represent the model’s tangent-linear operators at time $t_k$, $e$ is the innovation error, and $K$ is the Kalman gain. The time update step involves integrating two simultaneous ordinary differential equations to propagate the state vector and the covariance matrix. The Joseph stabilizing form is used to compute the filtered covariance to retain the positive-definiteness and the symmetry of the covariance matrix.
\begin{algorithm} [h]
\caption{CD-EKF Algorithm}
\label{Alg:CD-EKF}
\begin{algorithmic}[1]
    \State $\hat{x}_{0|0} \leftarrow \hat{x}_0$, $P_{0|0} \leftarrow P_0$
    \While {$k \ge 0$}
    
        \State $\dot{\hat{x}}_k(t) = f(\hat{x}_k(t))$,  \hfill \break
         $\dot{P}_k(t) = A_k(t) P_k(t) + P_k(t) A^T_k(t) + \sigma_k(t) \sigma^T_k(t)$, \hfill \break for $t \in [t_{k}, t_{k+1}[$, where $\hat{x}_k(t_k) = \hat{x}_{k|k}$, $P_k(t_k) = P_{k|k}$
                
        \State $\hat{x}_{k+1|k} \leftarrow \hat{x}_k(t_{k+1})$

        \State $P_{k+1|k} \leftarrow P_k(t_{k+1})$
        
        \State $\hat{y}_{k+1|k} \leftarrow h_{k+1}(\hat{x}_{k+1|k})$
        
        \State $R_{e,k+1} \leftarrow C_{k+1} P_{k+1|k} C_{k+1} + R_{k+1}$
        
        \State $K_{k+1} \leftarrow P_{k+1|k} C^T_{k+1} R^{-1}_{e,k+1}$
        
        \State $e_{k+1} \leftarrow y_{k+1} - \hat{y}_{k+1|k}$

        \State $\hat{x}_{k+1|k+1} \leftarrow \hat{x}_{k+1|k} + K_{k+1} e_{k+1}$

        \State $ J_{k+1} \leftarrow I - K_{k+1} C_{k+1}$
        
        \State $P_{k+1|k+1} \leftarrow  J_{k+1} P_{k+1|k} J_{k+1}^T + K_{k+1} R_{e,k+1} K^T_{k+1}$
        
        \State $k \leftarrow k+1$

    \EndWhile
\end{algorithmic}
\end{algorithm}

\subsection{The Unscented Kalman Filter}
In the CD-UKF, an unscented transformation replaces the linearization step in the CD-EKF to form a Gaussian approximation by approximating the solution of the moment equations for the mean and covariance \citep{UT}. In the CD-UKF, $2(n_{\hat{x}}+n_w) + 1$ $\sigma$-points are deterministically-drawn such that they capture the mean and covariance of the original distribution. Then, the time update step is implemented by numerically integrating the drawn $\sigma$-points forward in time using the nonlinear dynamics. The Gaussian process approximation using this nonlinear transformation enables to capture higher-order moments better than the linearization-based Gaussian approximation in the CD-EKF. However, the $\sigma$-point required number of model simulations, which is $2(n_{\hat{x}}+n_w) + 1$ in the augmented CD system, makes the CD-UKF intractable in high-dimensional systems. Algorithm \ref{Alg:CD-UKF} shows the implementation steps for the time and measurement updates in the CD-UKF, where, $\alpha$, $\beta$, and $\kappa$ are tuning parameters that control the spread of the sigma points and the distribution of the state vector. $W_m$ and $W_c$ are the mean and covariance weights, respectively. $\mathcal{X}$ and $\mathcal{Y}$ are the computed sigma points for the augmented state vector and the predicted measurements, respectively. $\Bar{\mathcal{X}}$ contains the sigma points for the augmented state vector in the time update step such that, $\Bar{\mathcal{X}} \in \mathbb{R}^{n_{\hat{x}} \times (2n_{\hat{x}}+2n_w + 1)}$. $d\omega^{(i)}$ are the sigma points for the noise term. This is required in the time update step in order to propagate the samples through the stochastic dynamics. $\Bar{\sigma}\text{-}Points$ and $\sigma\text{-}Points$ are the sigma points computation routines in the time and measurement update steps, respectively. The sigma points are computed as described by \cite{Nielsen,NielsenCPC}.
\begin{algorithm}
\caption{CD-UKF Algorithm}
\label{Alg:CD-UKF}
\begin{algorithmic}[1]
    \State Given $\hat{x}_{0|0} \leftarrow \hat{x}_0$, $P_{0|0} \leftarrow P_0$, $\alpha$, $\beta$, $\kappa$, $\Bar{n} = n_{\hat{x}} + n_w$
    
    \State $\Bar{c} \leftarrow \alpha^2 (\Bar{n} + \kappa)$, $\Bar{\lambda} \leftarrow \alpha^2 (\Bar{n} + \kappa) - \Bar{n}$
    
    \State $\Bar{W}_m^{(0)} \leftarrow \frac{\Bar{\lambda}}{\Bar{n} + \Bar{\lambda}}$, $\Bar{W}_c^{(0)} \leftarrow \frac{\Bar{\lambda}}{\Bar{n} + \Bar{\lambda}}+(1-\alpha^2+\beta)$, \newline
    $\Bar{W}_m^{(i)}, \Bar{W}_c^{(i)} \leftarrow \frac{1}{2(\Bar{n} + \Bar{\lambda})}$, for $i \in \{1, 2, ..., 2\Bar{n}\}$ 

    \State $c \leftarrow \alpha^2 (n_{\hat{x}} + \kappa)$, $\lambda \leftarrow \alpha^2 (n_{\hat{x}} + \kappa) - n_{\hat{x}}$
    
    \State $W_m^{(0)} \leftarrow \frac{\lambda}{n_{\hat{x}} + \lambda}$, $W_c^{(0)} \leftarrow \frac{\lambda}{n_{\hat{x}} + \lambda}+(1-\alpha^2+\beta)$, \newline
    $W_m^{(i)}, W_c^{(i)} \leftarrow \frac{1}{2(n_{\hat{x}} + \lambda)}$, for $i \in \{1, 2, ..., 2n_{\hat{x}}\}$ 

    \While {$k \ge 0$}
    
    \State $\Bar{\mathcal{X}}_{k|k}, d\omega^{(i)}_{k|k} \leftarrow \Bar{\sigma}\text{-}Points (\hat{x}_{k|k}, P_{k|k}, \Bar{c})$
    \For {$i = 1 : 2\Bar{n}+1$, and $t \in [t_k, t_{k+1}[$}
        \State $d{\Bar{\mathcal{X}}^{(i)}}_k(t) = f(\Bar{\mathcal{X}}^{(i)}_k(t)) dt + \sigma(\Bar{\mathcal{X}}^{(i)}_k(t)) d\omega^{(i)}_{k|k}(t)$, where $\Bar{\mathcal{X}}_k(t_k) = \Bar{\mathcal{X}}_{k|k}$
        
        \State $\Bar{\mathcal{X}}_{k+1|k} \leftarrow \Bar{\mathcal{X}}_k(t_{k+1})$

    \EndFor
        
        \State $\hat{x}_{k+1|k} \leftarrow \Bar{W}_m \Bar{\mathcal{X}}_{k+1|k}$
        
        \State $\Bar{E}^i_{x, k+1} \leftarrow \Bar{\mathcal{X}}^{(i)}_{k+1|k} - \hat{x}_{k+1|k}$ \Comment{the \textit{i}$^{th}$ column of $\Bar{E}_x$}

        \State $P_{k+1|k} \leftarrow \Bar{W}_c \Bar{E}_{x, k+1} \Bar{E}_{x, k+1}^T$

        \State $\mathcal{X}_{k+1|k} \leftarrow \sigma\text{-}Points (\hat{x}_{k+1|k}, P_{k+1|k}, c)$

        \State $\mathcal{Y}_{k+1|k} \leftarrow h_{k+1}(\mathcal{X}_{k+1|k})$
        
        \State $\hat{y}_{k+1|k} \leftarrow W_m \mathcal{Y}_{k+1|k}$

        \State $e_{k+1} \leftarrow y_{k+1} - \hat{y}_{k+1|k}$

        \State $E^i_{x, k+1} \leftarrow \mathcal{X}^{(i)}_{k+1|k} - \hat{x}_{k+1|k}$ \Comment{the \textit{i}$^{th}$ column of $E_x$}

        \State $E^i_{y, k+1} \leftarrow \mathcal{Y}^{(i)}_{k+1|k} - \hat{y}_{k+1|k}$ \Comment{the \textit{i}$^{th}$ column of $E_y$}

        \State $R_{xy, k+1} \leftarrow W_c E_{x, k+1} E_{y, k+1}^T$

        \State $R_{yy, k+1} \leftarrow W_c E_{y, k+1} E_{y, k+1}^T$
        
        \State $K_{k+1} \leftarrow R_{xy, k+1} R^{-1}_{yy, k+1}$

        \State $\hat{x}_{k+1|k+1} \leftarrow \hat{x}_{k+1|k} + K_{k+1} e_{k+1}$

        \State $P_{k+1|k+1} \leftarrow P_{k+1|k} - K_{k+1} R_{yy, k+1} K_{k+1}^T$

        \State $k \leftarrow k+1$
    
    \EndWhile
\end{algorithmic}
\end{algorithm}

\subsection{The Ensemble Kalman Filter}
The EnKF is a MC-sampling-based filtering approach that uses a finite number of randomly-drawn samples to represent the error statistics. The EnKF aims mainly to simplify the computational complexity in DA problems in high-dimensional nonlinear systems. Instead of propagating the full-rank covariance matrix, a low-rank sample covariance is approximated from the individually propagated samples. This allows for solving the estimation problem in an ensemble subspace instead of the original model space. However, the propagation of a low-rank covariance matrix can lead to spurious long-range correlations between uncorrelated state variables and causes an erroneous reduction in the estimated uncertainty, which can eventually lead to ensemble collapse and filter divergence. Techniques such as localization and covariance inflation can be implemented to overcome this problem. Algorithm \ref{Alg:CD-EnKF} illustrates the implementation steps for the EnKF, where $N$ is the number of ensemble members, $A'$ is the ensemble perturbation matrix, and $\mathbf{1}_N \in \mathbb{R}^{N \times N}$ is a matrix with all elements equal to $1/N$. 
\begin{algorithm}
\caption{CD-EnKF Algorithm}
\label{Alg:CD-EnKF}
\begin{algorithmic}[1]
    \State $\hat{x}^{(i)}_{0|0} \sim \mathcal{N}(\hat{x}_0, P_0), i = \{1, 2, .., N\}$
    \While {$k \ge 0$}
    
    \For {$i = 1 : N$, and $t \in [t_k, t_{k+1}[$}
        \State $d{\hat{x}}^{(i)}_k(t) = f(\hat{x}^{(i)}_k(t)) dt + \sigma(\hat{x}^{(i)}_k(t)) d\omega_k(t)$, where $\hat{x}^{(i)}_k(t_k) = \hat{x}^{(i)}_{k|k}$
        \State $\hat{x}^{(i)}_{k+1|k} \leftarrow \hat{x}^{(i)}_k(t_{k+1})$
    \EndFor
    
        \State $\hat{x}_{k+1|k} \leftarrow \frac{1}{N} \sum^N_{i=1} \hat{x}^{(i)}_{k+1|k}$

        \State $A_{k+1|k} \leftarrow [\hat{x}^{(1)}_{k+1|k}, \hat{x}^{(2)}_{k+1|k}, ..., \hat{x}^{(N)}_{k+1|k}]$

        \State $\Bar{A}_{k+1|k} \leftarrow A_{k+1|k} \mathbf{1}_N$

        \State $A_{k+1|k}' \leftarrow A_{k+1|k}- \Bar{A}_{k+1|k}$
        
        \State $P_{k+1|k} \leftarrow \frac{A_{k+1|k}' (A_{k+1|k}')^T}{N-1}$

    \For {$i = 1 : N$}

        \State $\hat{y}^{(i)}_{k+1|k} \leftarrow h_{k+1}(\hat{x}^{(i)}_{k+1|k})$

        \State $y^{(i)}_{k+1} \leftarrow y_{k+1} + v^{(i)}_{k+1}$
        
        \State $e^{(i)}_{k+1} \leftarrow y^{(i)}_{k+1} - \hat{y}^{(i)}_{k+1|k}$

    \EndFor
        \State $A_{y, k+1|k} \leftarrow [\hat{y}^{(1)}_{k+1|k}, \hat{y}^{(2)}_{k+1|k}, ..., \hat{y}^{(N)}_{k+1|k}]$

        \State $\Bar{A}_{y, k+1|k} \leftarrow A_{y, k+1|k} \mathbf{1}_N$

        \State $A_{y, k+1|k}' \leftarrow A_{y, k+1|k} - \Bar{A}_{y, k+1|k}$
        
        \State $R_{xy, k+1} \leftarrow \frac{A_{k+1|k}' (A_{y, k+1|k}')^T}{N-1}$

        \State $R_{yy, k+1} \leftarrow \frac{A_{y, k+1|k}' (A_{y, k+1|k}')^T}{N-1} + R_{k+1}$

        \State $K_{k+1} \leftarrow R_{xy, k+1} R^{-1}_{yy, k+1}$

    \For {$i = 1 : N$}
        \State $\hat{x}^{(i)}_{k+1|k+1} \leftarrow \hat{x}^{(i)}_{k+1|k} + K_{k+1} e^{(i)}_{k+1}$
    \EndFor
    
        \State $\hat{x}_{k+1|k+1} \leftarrow \frac{1}{N} \sum^N_{i=1} \hat{x}^{(i)}_{k+1|k+1}$

        \State $A_{k+1|k+1} \leftarrow [\hat{x}^{(1)}_{k+1|k+1}, \hat{x}^{(2)}_{k+1|k+1}, ..., \hat{x}^{(N)}_{k+1|k+1}]$

        \State $\Bar{A}_{k+1|k+1} \leftarrow A_{k+1|k+1} \mathbf{1}_N$

        \State $A_{k+1|k+1}' \leftarrow A_{k+1|k+1}- \Bar{A}_{k+1|k+1}$
        
        \State $P_{k+1|k+1} \leftarrow \frac{A_{k+1|k+1}' (A_{k+1|k+1}')^T}{N-1}$

        \State $k \leftarrow k+1$
    
    \EndWhile
\end{algorithmic}
\end{algorithm}

\subsection{The Particle Filter}
All the previous Kalman-based filters involve a form of Gaussian approximation, which inherently assumes linear system dynamics and Gaussian noises. The PF, on the other hand, is a fully nonlinear DA approach that does not postulate assumptions on the dynamic model or the state distribution. It only introduces one approximation to solving the FPE which is the MC sampling representation for the full uncertainty distribution, and therefore, the PF precisely samples the Bayesian posterior in the general case. Each MC particle in the PF, however, is assigned a weight that is updated to correct the difference between the approximated and the target distribution, resulting in the sequential importance resampling (SIR) procedure used in the PF \citep{sarkka_2013}. Similar to the EnKF, the PF can experience weight collapse, which can be mitigated using a resampling step. Algorithm \ref{Alg:CD-PF} indicates the implementation steps for the PF using SIR, where $w$ represents the likelihood weights.

Despite the perspicuous superiority of the PF over the filters that use Gaussian approximation due to its generalization properties, the Kalman-based filters are still the most commonly used in applications. This is due to the computational cost of the large number of model simulations required in the PF, which makes it infeasible, or at least computationally expensive, in a wide range of practical problems. In addition, most of the models in real applications can be fairly approximated as Gaussian processes. In some cases, however, the Kalman-based filters can even overperform the PF \citep{SimoSarrkaPhD}.
\begin{algorithm}
\caption{CD-PF Algorithm}
\label{Alg:CD-PF}
\begin{algorithmic}[1]
    \State $\hat{x}^{(i)}_{0|0} \leftarrow MC Sampling(p(\hat{x}(t_0))), i = \{1, 2, .., n_p\}, w_0$
    \While {$k \ge 0$}
    
    \For {$i = 1 : n_p$, and $t \in [t_k, t_{k+1}[$}
        \State $d{\hat{x}}^{(i)}_k(t) = f(\hat{x}^{(i)}_k(t)) dt + \sigma(\hat{x}^{(i)}_k(t)) d\omega_k(t)$, where $\hat{x}^{(i)}_k(t_k) = \hat{x}^{(i)}_{k|k}$
        \State $\hat{x}^{(i)}_{k+1|k} \leftarrow \hat{x}^{(i)}_k(t_{k+1})$
                
        \State $\hat{y}^{(i)}_{k+1|k} \leftarrow h_{k+1}(\hat{x}^{(i)}_{k+1|k})$
        
        \State $e^{(i)}_{k+1} \leftarrow y_{k+1} - \hat{y}^{(i)}_{k+1|k}$

        \State $w^{(i)}_{k+1} \propto w^{(i)}_{k} p(y_{k+1}|\hat{x}^{(i)}_{k+1})$
        
    \EndFor
    
    \For{$i = 1 : n_p$}
    
        \State $\Tilde{w}^{(i)}_{k+1} \leftarrow \frac{w^{(i)}_{k+1}}{\sum^{n_p}_{j=1} w^{(j)}_{k+1}}$ \Comment{Normalized weights}

    \EndFor
        
        \State $\hat{x}_{k+1|k+1} \leftarrow \sum^{n_p}_{i=1} \Tilde{w}^{(i)}_{k+1} \hat{x}^{(i)}_{k+1|k+1}$

        \State $E^i_{x, k+1} \leftarrow \hat{x}^{(i)}_{k+1|k} - \hat{x}_{k+1|k}$ \Comment{the \textit{i}$^{th}$ column of $E_x$}

        \State $P_{k+1|k+1} \leftarrow \frac{E_{x, k+1} E_{x, k+1}^T}{n_p-1}$ \Comment{Optional}

        \State $k \leftarrow k+1$
    
    \EndWhile
\end{algorithmic}
\end{algorithm}

\section{Numerical Case Study} \label{sec:App}
In this section, a second-order irreversible exothermic reaction conducted in an adiabatic CSTR is used to investigate the performance of the different CD assimilation methods in dual estimation problems.

\subsection{CSTR System}
The CSTR is a chemical reactor that is used for perfect mixing in liquid-phase or multiphase reactions. The nonlinear dynamics of the CSTR are described as follows
\begin{subequations}
\label{eq:odeCSTR}
\begin{align}
    \dot{C}_A &= \frac{F}{V} (C_{A, in} - C_A) - k(T) C_A C_B, \\
    \dot{C}_B &= \frac{F}{V} (C_{B, in} - C_B) - 2 k(T) C_A C_B, \\
    \dot{T} &= \frac{F}{V} (T_{in} - T) + \beta k(T) C_A C_B,
\end{align}
\end{subequations}
where $T$ is the temperature of the reactor, $C_A$ and $C_B$ are the concentration of substances $A$ and $B$ in the reactor, respectively. The rate constant $k(T) = k_0 exp(\frac{-E_a}{RT})$, and the parameter $\beta$ is defined from the reaction enthalpy as $\beta = -\Delta H_r/(\rho c_p)$.  Table \ref{tbl:CSTR} defines the remaining model variables and parameters along with their values.

The dynamics of the CSTR can be modeled using a SDE to capture random variations in the input temperature \citep{Jorgensen}. The resulting model will be the same as given in \eqref{eq:odeCSTR} for the first two states, $C_A$ and $C_B$; however, the temperature dynamic equation will include a stochastic term and will be represented in the form of \eqref{eq:CDStochModel} as follows
\begin{equation}
    dT = \left(\frac{F}{V} (T_{in} - T) + \beta k(T) C_A C_B\right) dt + \frac{F}{V} \sigma_T d\omega(t).
    \label{eq:sdeCSTR}
\end{equation}

\begin{figure}[t]
\begin{center}
\includegraphics[width=8.4cm]{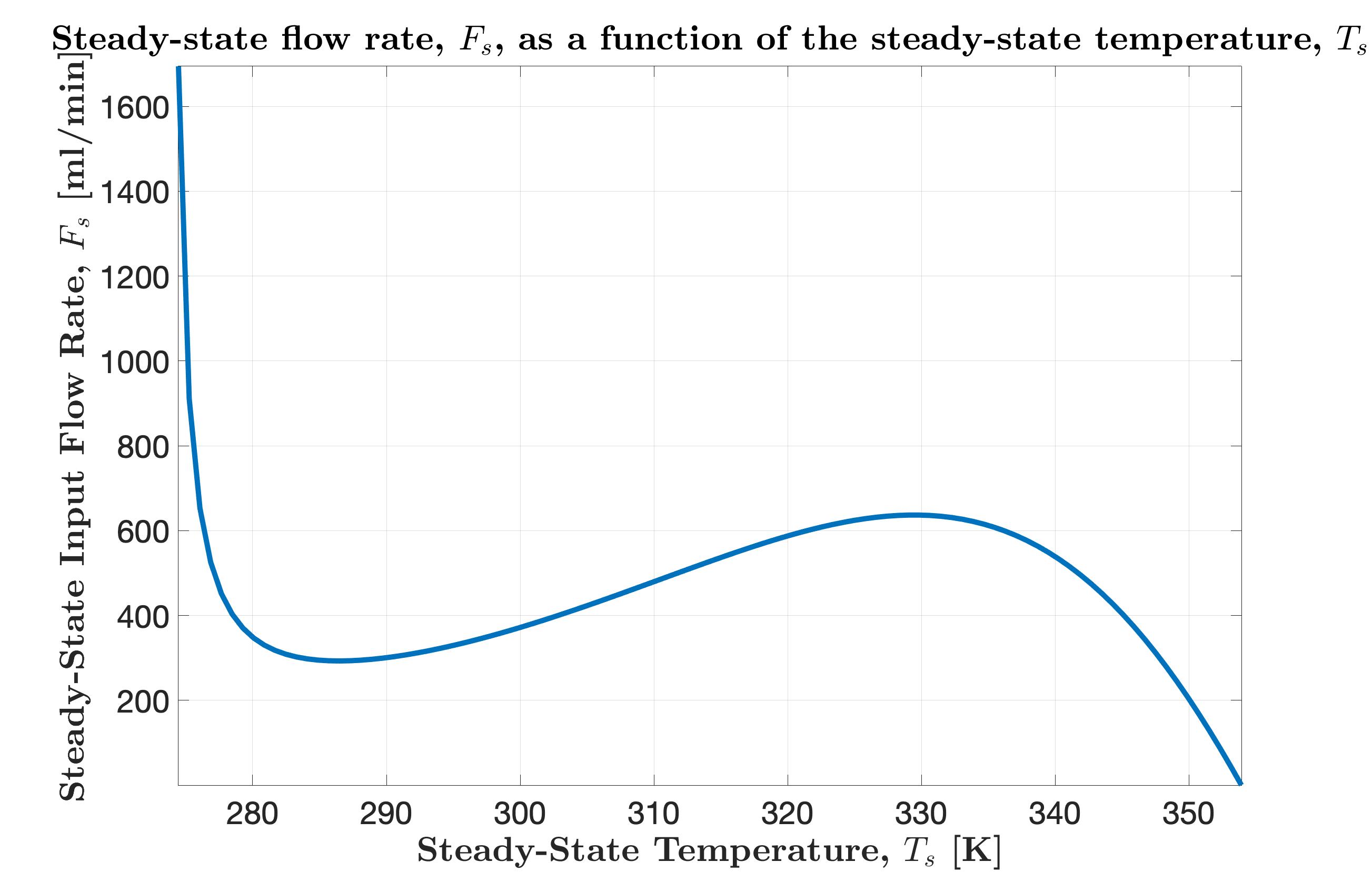}
\caption{Steady-state flow rate, $F_s$, as a function of the steady-state temperature, $T_s$.} 
\label{fig:ssFT}
\end{center}
\end{figure}

\begin{figure}[t]
\begin{center}
\includegraphics[width=8.4cm]{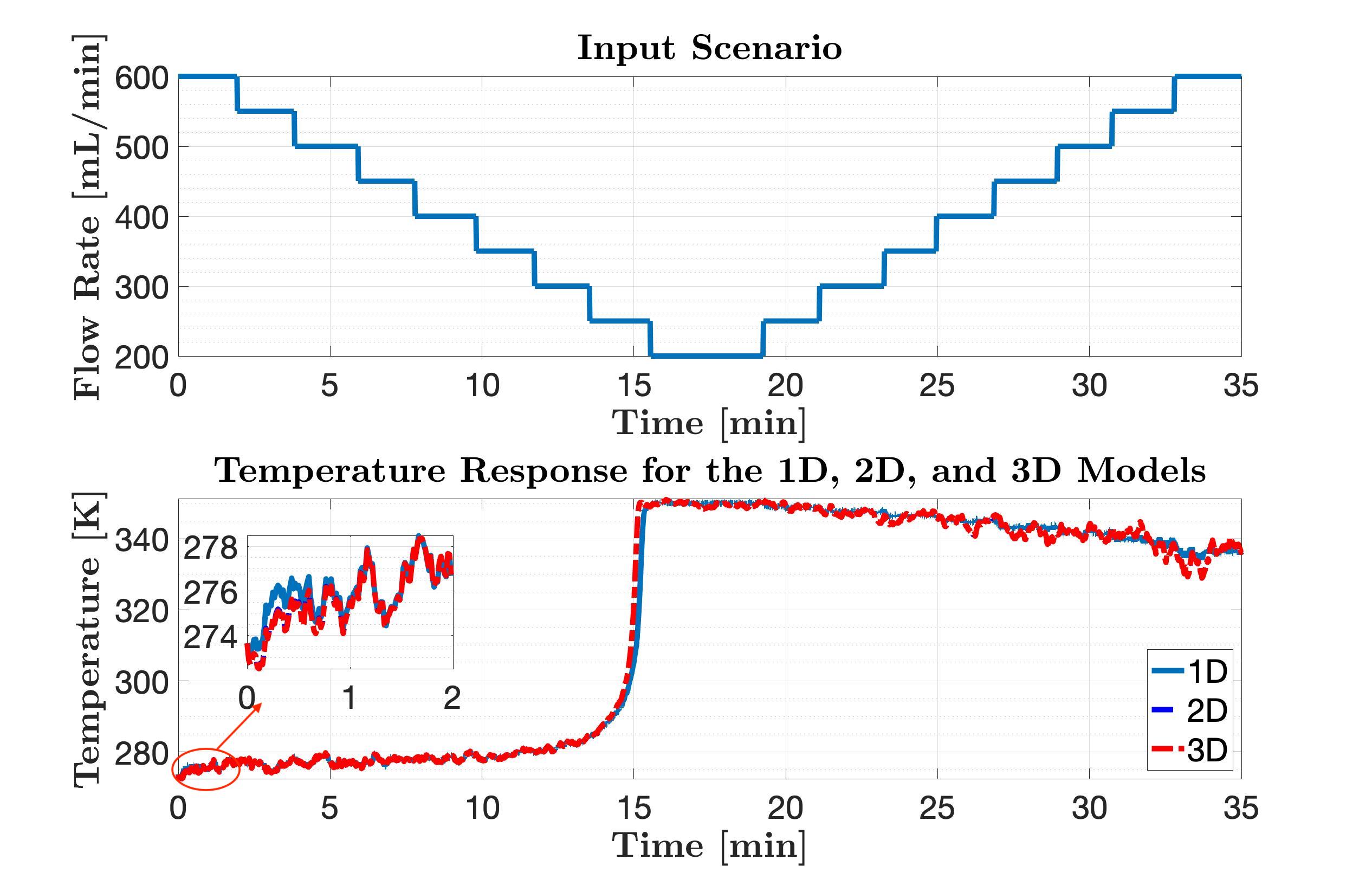}
\caption{Temperature response for the three-, two-, and one-dimensional representation of the stochastic CSTR model to a given input trajectory, $F$.} 
\label{fig:T_1d2d3d}
\end{center}
\end{figure}

The CSTR system is an interesting benchmark in the control and estimation field due to its interesting characteristics. For example, the CSTR model has multiple stable and unstable steady states in its operating window \citep{WAHLGREEN}. This is illustrated in Fig.~\ref{fig:ssFT} where the steady-state flow rate, $F_s$, is plotted as a function of the steady-state temperature, $T_s$. Moreover, using the definition of the extent of reaction, $X(t) = \frac{C_{B, in} - C_B(t)}{C_{B, in}}$, the three-state model, ($C_A(t)$, $C_B(t)$, and $T(t)$), can be approximated using two- and one-state reduced-order models in terms of ($X(t)$, and $T(t)$), or ($T(t)$), respectively \citep{WAHLGREEN}. Fig.~\ref{fig:T_1d2d3d} shows one realization of the simulated temperature response of the stochastic CSTR system for a given input trajectory, $F$, using the three\nobreakdash-, two-, and one-dimensional models.

\subsection{Estimation Twin Experiment}
The four CD filters, CD-EKF, CD-UKF, CD-EnKF, and CD-PF, are implemented to estimate the states, $C_A(t)$, $C_B(t)$, and $T(t)$, and the parameters, $\beta$ in this case, of the three-dimensional CSTR stochastic model. An initial state vector $x_0 = [0, 0, 273.65]^T$ and a parameter $\beta = 133.7792$ were selected as the ground truth model. The diffusion term $\sigma_T$ in \eqref{eq:sdeCSTR} is set equal to $5$. The input trajectory of the flow rate, $F$, is generated as shown in Fig. \ref{fig:T_1d2d3d}. The remaining parameters are defined as given in Table~\ref{tbl:CSTR}. Synthetic measurements for the temperature are simulated using the true model by adding measurement noise $v(t_k) \sim \mathcal{N}_{iid}(0, R_k)$, with $R_k = \sigma^2_v = 3^2$. The Euler-Maruyama integration scheme is used to simulate the stochastic system. The model was simulated for $35$ minutes with 210 equidistant samples. Then, the four filters are implemented to estimate the three states and the parameter by assimilating the generated noisy observations of the temperature. The initial guess for the state vector was $x_{e0} = [0.1, 0.2, 293.65]^T$, and for the parameter was $\beta_{e0} = 123.7792$. The tuning parameters for the UKF are $[\alpha, \beta, \kappa] = [0.2, 2, 0]$. For both the EnKF and the PF, the ensemble size is $1000$ MC samples, initiated by drawing from an initial multivariate Gaussian distribution for the augmented state vector.

\begin{table}[t]
    \caption{CSTR Model Parameters}
    \label{tbl:CSTR}
\begin{adjustbox}{width=\columnwidth,center}
        \begin{tabular}{l l l l}
            \hline
            Symbol                  & Variable                                    & Value     & Unit \\
            \hline
            $\rho$                  & Density                       & $1.0$      & $kg/L$ \\
            $c_p$                    & Specific heat capacity    & $4.186$        & $kJ/(kg.K)$ \\
            $k_0$                   & Arrhenius constant               & $exp(24.6)$ &$L/(mol\text{·}s)$ \\
            $E_a/R$                 & Activation energy                   & $8500$ & $K$ \\
            $\Delta H_r$            & Reaction enthalpy      & $9.0 \times 10^{-4}$ & Pa.s \\
            $V$                     & Reactor volume                   & $0.105$ &$L$ \\
            $C_{A, in}$             & Feed concentration, A           & $1.6/2$ & $mol/L$ \\
            $C_{B, in}$             & Feed concentration, B           & $2.4/2$ & $mol/L$ \\
            $T_{in}$               & Inlet temperature              & $273.65$ &$K$ \\
            \hline
        \end{tabular}
    \end{adjustbox}
\end{table}

\begin{figure} [b]
\begin{center}
\includegraphics[width=8.4cm]{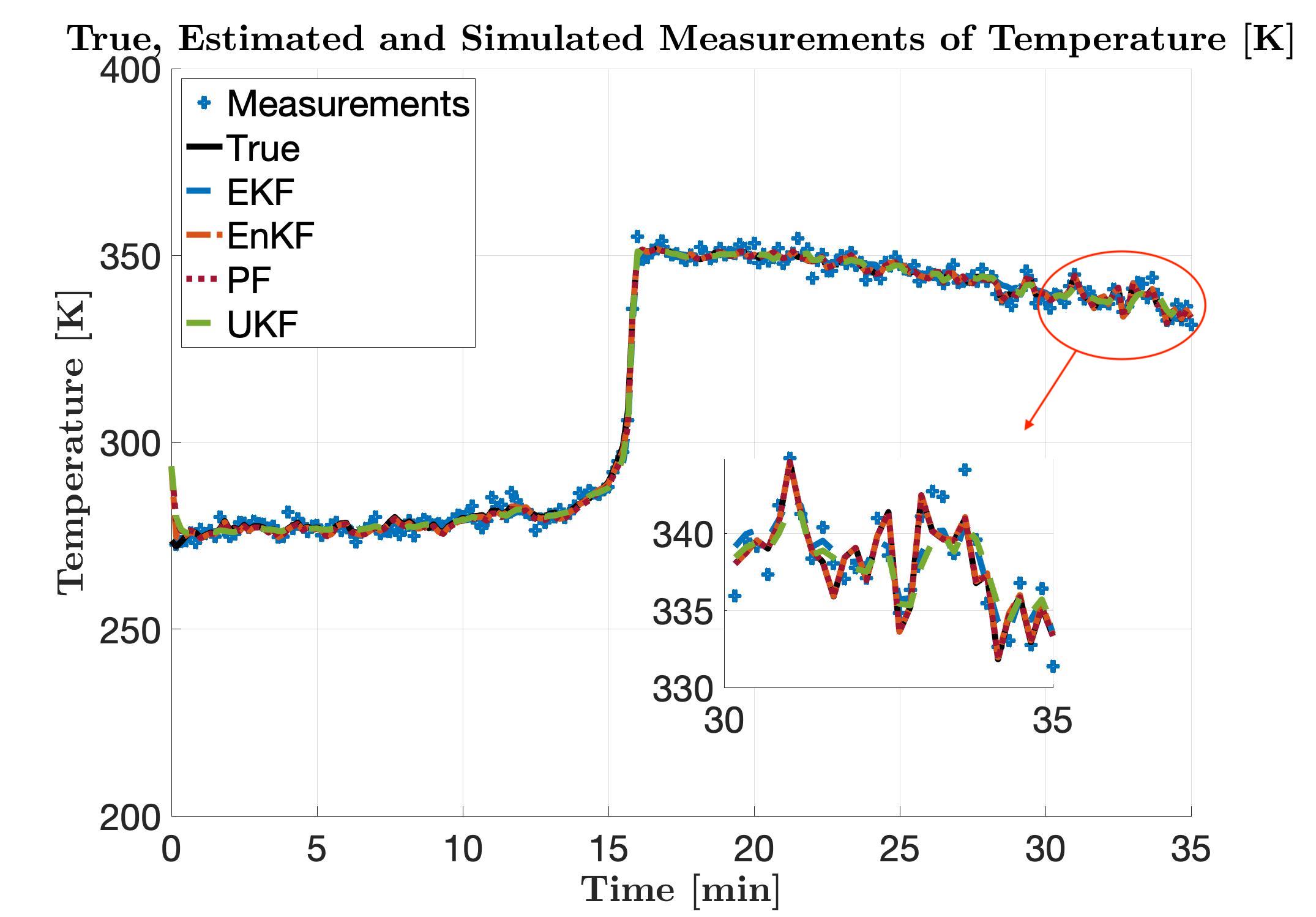}
\caption{True and estimated trajectories along with the measurements for the temperature, $T$.}
\label{fig:AllFilters_T}
\end{center}
\end{figure}

\begin{figure}[b]
\begin{center}
\includegraphics[width=8.4cm]{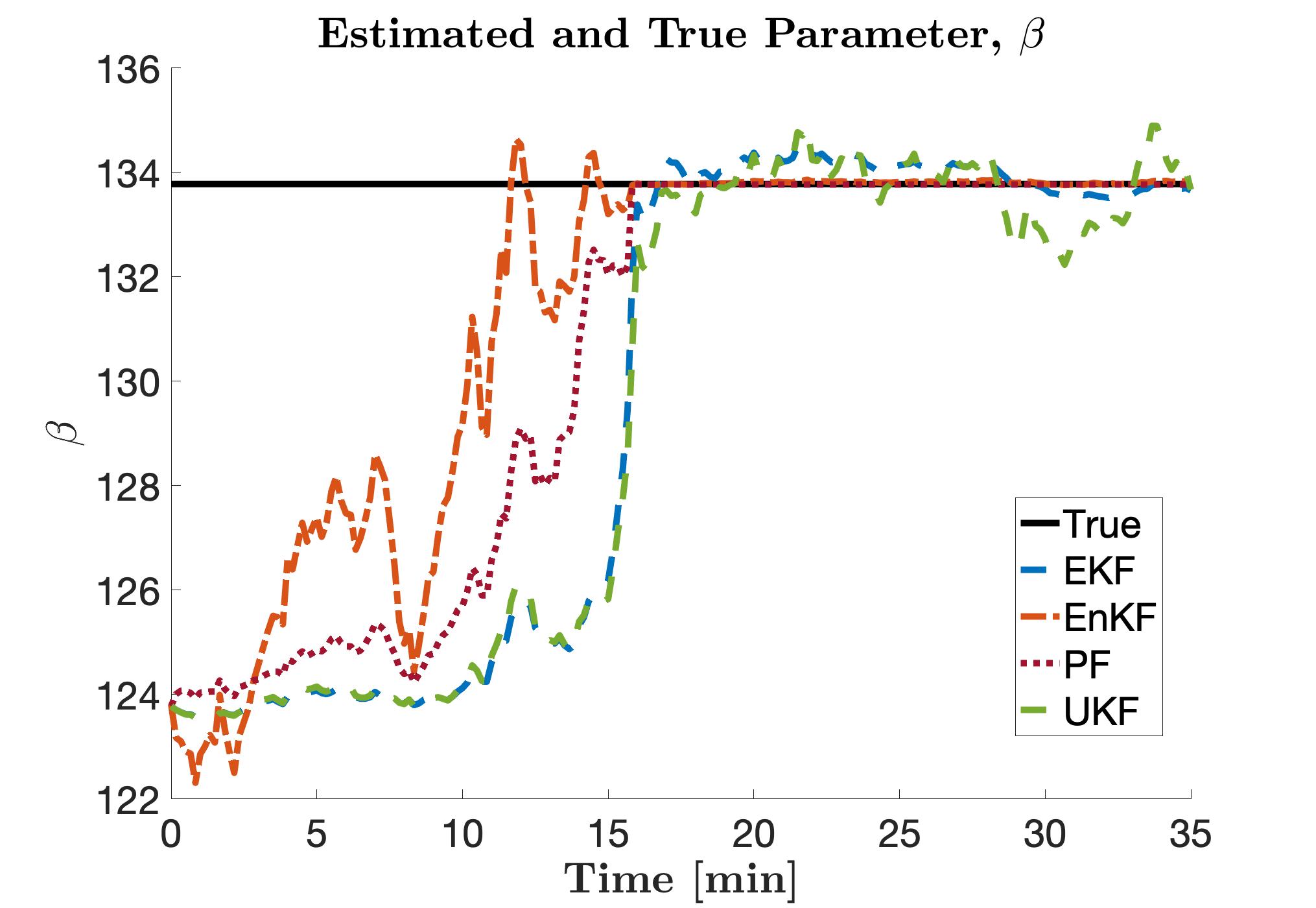}
\caption{True and estimated trajectories for the parameter, $\beta$.}
\label{fig:AllFilters_Beta}
\end{center}
\end{figure}

\subsection{Results and Discussion}
\begin{figure}[b]
\begin{center}
\includegraphics[width=8.4cm]{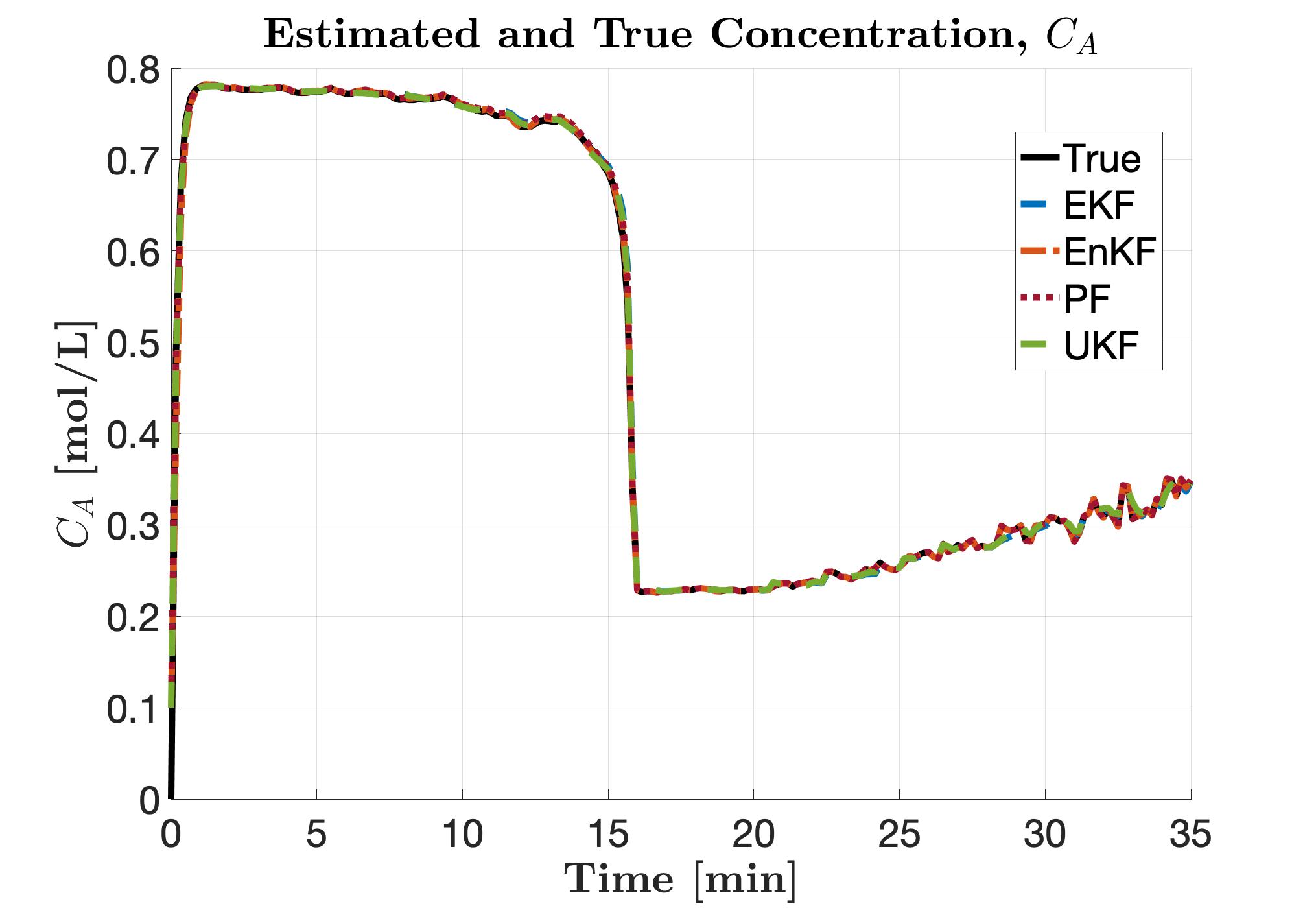}
\caption{True and estimated trajectories for the concentration of $A$, $C_A$.}
\label{fig:AllFilters_Ca}
\end{center}
\end{figure}

\begin{figure} [b]
\begin{center}
\includegraphics[width=0.44\textwidth]{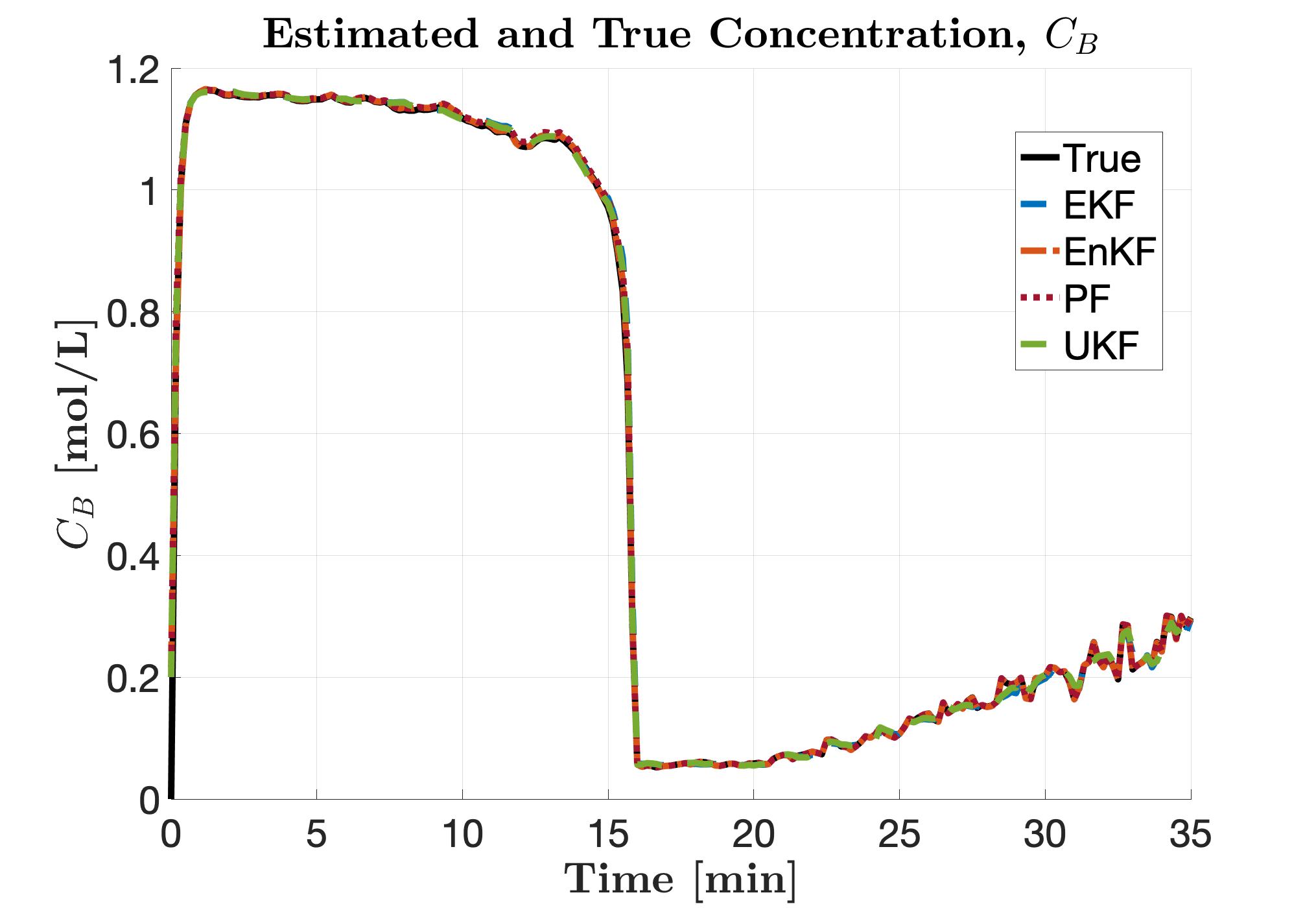}
\caption{True and estimated trajectories for the concentration of $B$, $C_B$.}
\label{fig:AllFilters_Cb}
\end{center}
\end{figure}

Fig.~\ref{fig:AllFilters_T} -~\ref{fig:AllFilters_Cb} illustrate the assimilation results for all the filters in estimating the four unknown variables. The results show that the four filters can track the state trajectories fairly well, with the EKF and the UKF showing some small discrepancies at some points; while the EnKF and the PF, on the other hand, exhibit nearly zero-error tracking in the whole assimilation window. In the parameter estimation problem, however, the estimation starts from an initial guess that is far away from the true value for the parameter $\beta$. All the filters then start to converge gradually until the estimation error is eliminated rapidly around the point when the system states exhibit rapid changes. Similarly, it can be seen that the estimates of the EnKF and the PF remain steady at the true parameter value with no more updates after that point, while the EKF and the UKF keep fluctuating around the true value.

\subsubsection{Estimation Errors}
Table~\ref{tbl:Errors} summarizes the estimation performance of the four filters quantitatively. It shows the MSE for the states, $MSE_x$, for the parameter, $MSE_x$, and the computational cost, $t_{CPU}$, for each filter. It indicates that the EnKF, in this case, gives the lowest estimation error for both the state and the parameter estimation problems. The EKF and the UKF perform equally well in terms of the MSE in this case, while the MC-based filters result in lower estimation errors.

\subsubsection{Computational Cost}
Table~\ref{tbl:Errors} presents the CPU computation time, $t_{CPU}$, for a single assimilation step for each estimation approach. It indicates that the most computationally efficient approach, in this case, is the EKF, which is a reasonable result in this low-dimensional case with precomputed analytical Jacobians. Compared to other $\sigma$-point and ensemble- or particle-based filters, which require a number of model simulations at each assimilation step, the EKF performs simpler steps to propagate the covariance information. However, this might not be the case in high-dimensional systems when it becomes computationally expensive to propagate the full-rank covariance matrix, or when the analytical Jacobians are not available and have to be computed numerically at each assimilation step. In addition, the computation times show that the MC sampling approaches are considerably more computationally expensive due to the large number of required simulations at each assimilation step.

\subsubsection{Uncertainty Quantification}
As a fully nonlinear DA approach, a CD-PF with a large number of particles ($10000$) was used to approximate the correct uncertainty propagation in order to investigate the performance of each filter in quantifying the uncertainty. The results show that the four filters can approximate the correct covariance in this case with arbitrarily small errors. However, one interesting observation on the ensemble size in the MC-based filters is that the performance of the EnKF is more robust than the PF to changes in the ensemble size. That is, decreasing the number of the ensemble members in the EnKF to a considerably smaller number does not degrade its performance. This was the case until an ensemble size of the same dimensionality as the state vector, which is $4$ in this example. While the performance of the PF, on the other hand, degrades drastically when the number of the particles was decreased to $100$. This can be attributed to the inherently nonlinear propagation of the uncertainty in the PF, which can be more sensitive to the number of the MC samples. This also indicates that the EnKF results in an equivalent estimation performance to the PF, or even better in this case, with a much smaller number of MC simulations, and therefore, with better computational efficiency. The EnKF experiences an ensemble collapse and a filter divergence when the ensemble size becomes smaller than the number of variables in the augmented state vector.

\begin{table}
    \caption{Estimation Performance for the Four Filters}
    \label{tbl:Errors}
    \begin{center}
        \begin{tabular}{l l l l l}
            \hline
            Index   & EKF       & UKF       & EnKF       & PF \\
            \hline
            $MSE_x$             & 0.6628    & 0.6377    & 0.4795     & 0.5251 \\
            $MSE_p$             & 6.2995    & 6.3048    & 4.5204     & 4.9463 \\
            $t_{CPU} [s]$            & 0.0297    & 0.0821    & 4.0285     & 4.1665 \\
            \hline
        \end{tabular}
    \end{center}
\end{table}

\section{Conclusion} \label{sec:conc}
In this study, four continuous-discrete filters, the CD-EKF, CD-UKF, CD-EnKF, and CD-PF, are presented as approximate solutions to the FPE using different assumptions. The characteristics and the limitations of each filter are discussed and illustrated in a combined parameter and state estimation numerical example. A CSTR model described by a SDE is used to illustrate and compare the performance of the different approaches. The results show that the CD-EKF is the most computationally efficient filter in this case, while the CD-EnKF results in the smallest estimation errors.


\bibliography{ifacconf}             

\begin{thebibliography}{17}
\providecommand{\natexlab}[1]{#1}
\providecommand{\url}[1]{\texttt{#1}}
\providecommand{\urlprefix}{URL }
\expandafter\ifx\csname urlstyle\endcsname\relax
  \providecommand{\doi}[1]{doi:\discretionary{}{}{}#1}\else
  \providecommand{\doi}{doi:\discretionary{}{}{}\begingroup
  \urlstyle{rm}\Url}\fi

\bibitem[{Bui-Thanh(2021)}]{bui2021optimality}
Bui-Thanh, T. (2021).
\newblock The optimality of bayes' theorem.
\newblock \emph{SIAM news}, 54(6).

\bibitem[{Bui-Thanh and Ghattas(2015)}]{BayesOptimal}
Bui-Thanh, T. and Ghattas, O. (2015).
\newblock {B}ayes is optimal.
\newblock \emph{ICES REPORT, The Institute for Computational Engineering and
  Sciences, The University of Texas at Austin}, 15(04).

\bibitem[{Carrassi et~al.(2018)Carrassi, Bocquet, Bertino, and
  Evensen}]{carrassi2018data}
Carrassi, A., Bocquet, M., Bertino, L., and Evensen, G. (2018).
\newblock Data assimilation in the geosciences: An overview of methods, issues,
  and perspectives.
\newblock \emph{Wiley Interdisciplinary Reviews: Climate Change}, 9(5), e535.

\bibitem[{Diaa-Eldeen et~al.(2022)Diaa-Eldeen, Berg, and Hovd}]{Diaa-Eldeen}
Diaa-Eldeen, T., Berg, C.F., and Hovd, M. (2022).
\newblock Observability-aware ensemble {K}alman filter for reservoir model
  updating.
\newblock In \emph{2022 European Control Conference (ECC)}, 1714--1721.
\newblock \doi{10.23919/ECC55457.2022.9838071}.

\bibitem[{Evensen et~al.(2022)Evensen, Vossepoel, and van
  Leeuwen}]{evensen2022data}
Evensen, G., Vossepoel, F.C., and van Leeuwen, P.J. (2022).
\newblock \emph{Data {A}ssimilation {F}undamentals: A Unified Formulation of
  the State and Parameter Estimation Problem}.
\newblock Springer Nature.

\bibitem[{Gunther et~al.(1997)Gunther, Beard, Wilson, Oliphant, and
  Stirling}]{GalerkinMethod}
Gunther, J., Beard, R., Wilson, J., Oliphant, T., and Stirling, W. (1997).
\newblock Fast nonlinear filtering via {G}alerkin's method.
\newblock In \emph{Proceedings of the 1997 American Control Conference},
  volume~5, 2815--2819 vol.5.

\bibitem[{Julier and Uhlmann(2004)}]{UT}
Julier, S.J. and Uhlmann, J.K. (2004).
\newblock Unscented filtering and nonlinear estimation.
\newblock \emph{Proceedings of the IEEE}, 92(3), 401--422.

\bibitem[{Jørgensen et~al.(2020)Jørgensen, Ritschel, Boiroux,
  Schroll-Fleischer, Wahlgreen, Krogh~Nielsen, Wu, and Huusom}]{Jorgensen}
Jørgensen, J.B., Ritschel, T.K.S., Boiroux, D., Schroll-Fleischer, E.,
  Wahlgreen, M.R., Krogh~Nielsen, M., Wu, H., and Huusom, J.K. (2020).
\newblock Simulation of {NMPC} for a laboratory adiabatic {CSTR} with an
  exothermic reaction.
\newblock In \emph{2020 European Control Conference (ECC)}, 202--207.

\bibitem[{Nielsen et~al.(2022)Nielsen, Ritschel, Christensen, Dragheim, Huusom,
  Gernaey, and Jørgensen}]{Nielsen}
Nielsen, M.K., Ritschel, T.K.S., Christensen, I., Dragheim, J., Huusom, J.K.,
  Gernaey, K.V., and Jørgensen, J.B. (2022).
\newblock State estimation methods for continuous-discrete nonlinear systems
  involving stochastic differential equations.

\bibitem[{Nielsen et~al.(2023)Nielsen, Ritschel, Christensen, Dragheim, Huusom,
  Gernaey, and Jørgensen}]{NielsenCPC}
Nielsen, M.K., Ritschel, T.K.S., Christensen, I., Dragheim, J., Huusom, J.K.,
  Gernaey, K.V., and Jørgensen, J.B. (2023).
\newblock State estimation for continuous-discrete-time nonlinear stochastic
  systems.
\newblock In \emph{Accepted in Foundations of Computer Aided Process Operations
  / Chemical Process Control Conference 2023 (FOCAPO/CPC), Texas, USA}.

\bibitem[{Pavliotis(2014)}]{Pavliotis:2014}
Pavliotis, G.A. (2014).
\newblock \emph{Stochastic Processes and Applications. Diffusion Processes, the
  Fokker-Planck and Langevin Equations}.
\newblock Springer, Heidelberg, Germany.

\bibitem[{Särkkä(2006)}]{SimoSarrkaPhD}
Särkkä, S. (2006).
\newblock \emph{{Recursive Bayesian inference on stochastic differential
  equations}}.
\newblock Doctoral thesis, Helsinki University of Technology.

\bibitem[{Särkkä(2013)}]{sarkka_2013}
Särkkä, S. (2013).
\newblock \emph{{B}ayesian Filtering and Smoothing}.
\newblock Institute of Mathematical Statistics Textbooks. Cambridge University
  Press.

\bibitem[{Thygesen(2022)}]{Thygesen:2022}
Thygesen, U.H. (2022).
\newblock \emph{Stochastic Differential Equations for Science and Engineering}.
\newblock Technical University of Denmark, Lecture Notes, Kgs. Lyngby, Denmark.

\bibitem[{Wahlgreen et~al.(2020)Wahlgreen, Schroll-Fleischer, Boiroux,
  Ritschel, Wu, Huusom, and J{\o}rgensen}]{WAHLGREEN}
Wahlgreen, M.R., Schroll-Fleischer, E., Boiroux, D., Ritschel, T.K., Wu, H.,
  Huusom, J.K., and J{\o}rgensen, J.B. (2020).
\newblock Nonlinear model predictive control for an exothermic reaction in an
  adiabatic {CSTR}.
\newblock \emph{IFAC-PapersOnLine}, 53(1), 500--505.

\bibitem[{Xu et~al.(2020)Xu, Zhang, Li, Zhou, Liu, and Kurths}]{xu2020solving}
Xu, Y., Zhang, H., Li, Y., Zhou, K., Liu, Q., and Kurths, J. (2020).
\newblock Solving {F}okker-{P}lanck equation using deep learning.
\newblock \emph{Chaos: An Interdisciplinary Journal of Nonlinear Science},
  30(1), 013133.

\bibitem[{Zellner(1988)}]{OptimalInformation}
Zellner, A. (1988).
\newblock Optimal information processing and {B}ayes's theorem.
\newblock \emph{The American Statistician}, 42(4), 278--280.

\end{thebibliography}
                                                   







\appendix
\end{document}